
\input epsf.tex
\input amssym.def
\input amssym
\magnification=1100
\baselineskip = 0.22truein
\lineskiplimit = 0.01truein
\lineskip = 0.01truein
\vsize = 8.7truein
\voffset = 0.2truein
\parskip = 0.10truein
\parindent = 0.3truein
\settabs 12 \columns
\hsize = 6.0truein
\hoffset = 0.1truein

\setbox\strutbox=\hbox{%
\vrule height .708\baselineskip
depth .292\baselineskip
width 0pt}
\font\caps=cmcsc10
\font\bigtenrm=cmr10 at 14pt

\def\sqr#1#2{{\vcenter{\vbox{\hrule height.#2pt
\hbox{\vrule width.#2pt height#1pt \kern#1pt
\vrule width.#2pt}
\hrule height.#2pt}}}}
\def\square{\mathchoice\sqr46\sqr46\sqr{3.1}6\sqr{2.3}4}

\centerline{\bigtenrm CORE CURVES OF TRIANGULATED SOLID TORI}
\tenrm
\vskip 14pt
\centerline{MARC LACKENBY}
\vskip 18pt
\centerline{\caps 1. Introduction}
\vskip 6pt

One of the most powerful tools in 3-manifold topology is normal
surface theory. Many topologically relevant surfaces can be placed
into normal form (or some variant of this) with respect to any
triangulation $T$ of the manifold $M$. An important
case of this phenomenon is when $M$ is the solid torus and
the surface is a meridian disc. The existence of a meridian disc
in normal form and the fact that this can be algorithmically
detected is the key to Haken's solution to the problem of recognizing
the unknot [1]. However, normal surface theory suffers from some 
substantial limitations, possibly the most important of which
is that every normal meridian disc may have `exponential complexity'.
More precisely, if $n$ is the number of tetrahedra in $T$,
then normal surface theory only produces a meridian disc $D$
with at most $2^{kn}$ normal triangles and squares, where $k = 10^{11}$
(see [2]). 
This is more than just an artifice of the theory, because one can find triangulations 
of the solid torus where every normal meridian disc has exponentially many squares
and triangles. Indeed, we will do this explicitly in Section 6.
However, in this paper, we show that 
for a related problem, there is a solution with linearly bounded complexity.
In addition to the meridian disc,
the solid torus contains another important sub-object:
the core curve $C$, which is defined, up to ambient isotopy,
to be $\{ \ast \} \times S^1 \subset D^2 \times S^1$, where $\ast$ is 
a point in the interior of $D^2$. In this paper,
we will address the problem of placing $C$ into `normal form'.
The surprising conclusion is that this can be
achieved with a linear upper bound on the `complexity' of
$C$. 

In order to state our main result most cleanly, we say that
a curve $C$ embedded in the solid torus $M$ is a {\sl pre-core curve}
if $C$ is a core curve in the manifold $M \cup (\partial M \times I)$
that is obtained from $M$ by attaching a collar to $\partial M$. Thus,
a pre-core curve need not lie in the interior of $M$.

Whenever we refer to a triangulation of a 3-manifold $M$, we use the more
general definition that is now standard in low-dimensional topology.
Thus, it is an expression of $M$ as a collection of 3-simplices,
with some of their faces identified in pairs, via affine homeomorphisms.
Since the gluing maps are affine, it makes
sense to speak of a straight line in a face of the triangulation.

The following is our main result.

\noindent {\bf Theorem 1.1.} {\sl Let $T$ be a triangulation
of the solid torus $M$. Then $M$ contains a pre-core curve $C$ that lies in
the 2-skeleton of $T$, intersects the edges in finitely many points, is disjoint
from the vertices and intersects the interior of each face in at most $10$ straight arcs.}

Note that this does not bound the number of intersections points between
$C$ and each edge, face or tetrahedron of $T$. This is because
as $C$ runs from one face to another, the result is a point of intersection
with an edge of $T$, and this may give rise to isolated points
of intersection with the faces and tetrahedra incident to that edge.
If one wants to bound the total number of points of intersection
between a core curve and each tetrahedron, it seems to be best to
make the curve transverse to the 2-skeleton of $T$, as follows.

\noindent {\bf Theorem 1.2.} {\sl Let $T$ be a triangulation
of the solid torus $M$. Then $M$ contains a core curve $C'$ that
intersects each tetrahedron $\Delta$ of $T$ in at most $18$ arcs. Moreover,
each such arc is properly embedded in $\Delta$, and has endpoints
in the interior of the faces of $\Delta$. The intersection $\Delta \cap C'$
is a trivial tangle in $\Delta$. In fact, $\Delta \cap C'$ is
parallel to a collection of arcs $\alpha$ in $\partial \Delta$, with
the property that the intersection between each component of $\alpha$
and each face of $\Delta$ is at most one straight arc.}

Thus, the above result asserts that there is a core curve that intersects
each tetrahedron in one of finitely many possible trivial tangles, where this
finite list is universal, in the sense that it is independent of the triangulation
$T$.

In addition, we will prove similar theorems for partially ideal triangulations
of the solid torus, and for `affine' handle structures. 
(See Theorems 4.1 and 4.2.)

Although Theorem 1.1 is a result about triangulations, it has geometric
consequences. The reason is that a triangulation of a manifold $M$ determines
a path metric in which each tetrahedron is regular and Euclidean,
and, conversely, any Riemannian metric can be approximated (in a suitable sense)
by such a piecewise Euclidean metric. By using this relationship between
triangulations and Riemannian metrics, and applying Theorem 1.1,
we obtain the following result.

\noindent {\bf Theorem 1.3.} {\sl For each $K$, $I > 0$, there is a constant
$c(K,I)$ with the following property. If $M$ is a solid torus with a Riemannian
metric having volume at most $V$, injectivity radius at least $I$ and all sectional
curvatures in the interval $(-K,K)$, then $M$ contains a core curve with length
at most $c(K,I) \ V$.}

It seems likely that Theorems 1.1 and 1.2 will have many other applications.
In fact, the version of these theorems for affine handle structures is a key step in the proof of the main theorem in [7].
This asserts that the crossing number of a satellite knot is at least
$10^{-13}$ times the crossing number of its companion knot. For more
details and an explanation of this terminology, see [7].

We now give an outline of the proof of Theorems 1.1 and 1.2.
We may find a meridian disc $D$ for the solid torus $M$ that is normal with respect to the
triangulation $T$. If we cut $M$ along $D$, the result is a 3-ball $X$.
The parts of $D$ lying between adjacent parallel normal discs patch up
to form an $I$-bundle in $X$, which we term the {\sl parallelity bundle} ${\cal B}$. 
(A precise definition is given in Section 2.)
Now, in each 3-simplex $\Delta$ of $T$, at most 6 components of $\Delta - N(D)$
do not lie in ${\cal B}$. The goal, therefore, is to find a core curve $C$ that
avoids ${\cal B}$ and that intersects each component of $\Delta - N(D)$ in a bounded
number of arcs with controlled topology. However, this is not completely straightforward.
For example, ${\cal B}$ may be almost all of $X$,
and so it is not clear that there is even a single core
curve disjoint from ${\cal B}$, particularly one with
controlled intersection with each tetrahedron.
The technique that we use is to construct a product structure $D^2 \times [-1,1]$ on
$X$, where $D^2 \times \{ -1, 1 \}$ is the copies of $D$ in $X$,
and with the property that the product structure agrees with the $I$-bundle structure on (most components
of) ${\cal B}$. This product structure determines a homeomorphism $\phi \colon D^2 \times \{1\}
\rightarrow D^2 \times \{ -1 \}$. We also have a homeomorphism $\psi \colon D^2 \times \{ -1 \}
\rightarrow D^2 \times \{ 1\}$ arising from the gluing map. Their composition $\psi \phi$ is
a homeomorphism $D^2 \times \{ 1 \} \rightarrow D^2 \times \{ 1 \}$. By the
Brouwer fixed point theorem, this has a fixed point $x \in D^2 \times \{ 1 \}$. There are
essentially two cases: either $x$ lies in ${\cal B}$ or it does not. If it
does, then $x \times [-1,1]$ is a vertical curve in ${\cal B}$, which
therefore joins adjacent parallel normal discs of $D$. But the two endpoints
of $x \times [-1,1]$ are identified under the gluing map, because $x$ is a fixed point of
$\psi \phi$, which is a contradiction.
Thus, $x$ does not lie in ${\cal B}$, and with some more work, one can arrange
that (a copy of) $x \times [-1,1]$ avoids ${\cal B}$. The endpoints of this arc patch together
to form the required core curve of $M$.

\vskip 18pt
\centerline{\caps 2. Normal surfaces and parallelity bundles}
\vskip 6pt

In this section, we recall some basic normal surface theory,
and then go on to define parallelity bundles. 

Recall that a disc properly embedded in a 3-simplex $\Delta$ is
{\sl normal} if it is disjoint from the vertices,
it intersects each edge transversely in at most one point, but it is not
disjoint from the edges. It is then
a {\sl triangle} or {\sl square}, as shown in Figure 1. 
A properly embedded surface in a triangulated 3-manifold is {\sl normal}
if it intersects each 3-simplex of the triangulation in a collection of normal discs.

\vskip 12pt
\centerline{
\epsfxsize=1.8in
\epsfbox{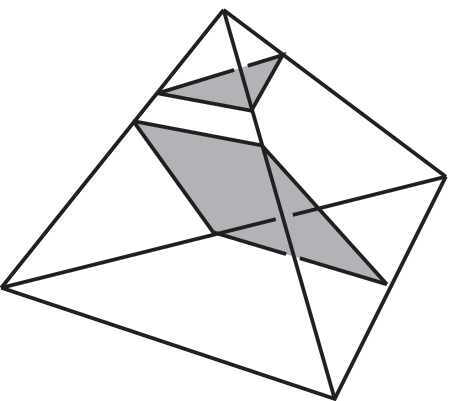}
}
\vskip 6pt
\centerline{Figure 1.}

Any normal surface may be straightened with respect to the triangulation, as follows.
First ambient isotope the arcs of intersection with the faces, so that they become straight.
Each triangle in each 3-simplex may then be ambient isotoped, keeping its boundary fixed,
so that it becomes planar in the affine structure. The same is not necessarily true of the squares.
But we may realise each square as the union of two flat pieces glued along a straight line.
We will henceforth assume that any normal surface has been straightened in this way.

An arc properly embedded in a 2-simplex is {\sl normal}
if it is disjoint from the vertices and has endpoints lying in
distinct edges. If $F$ is a compact surface with a triangulation,
then a 1-manifold properly embedded in $F$ is {\sl normal}
if it intersects each 2-simplex in a collection of normal arcs.
Define the {\sl length} of a normal 1-manifold in $F$
to be its number of intersections with the 1-skeleton.

The following is variant of a well-known fact in normal surface theory.
(See Proposition 4.4 in [6] for comparison.)

\noindent {\bf Proposition 2.1.} {\sl Let $T$ be a triangulation of
a compact irreducible 3-manifold $M$. Let $S$ be a properly embedded, oriented,
incompressible  surface in $M$, with no components that are 2-spheres
or boundary-parallel discs.
Give $\partial S$ the orientation that it inherits from $S$.
Suppose that $\partial S$ does not intersect any edge of $T \cap \partial M$ in 
two points with opposite orientations. Then, there is an ambient isotopy, 
keeping $\partial S$ fixed, taking $S$ into normal form.}

\noindent {\sl Proof.} In the usual normalisation procedure, there is
one step where $\partial S$ may need to be moved. This happens as
follows. Suppose that there is an arc of intersection between $S$
and a face of $T$, with endpoints in the same edge of the face,
and suppose that this edge lies in $\partial M$. The standard way to argue is to
assume that $S$ is boundary-incompressible and that $\partial M$
is incompressible, and thence perform an isotopy which removes
the arc. However, we do not follow this line of argument here.
Instead, we note that at the endpoints of the arc, $\partial S$
intersects the edge with opposite signs, which is contrary to
assumption. Thus, $\partial S$ does not need to be moved. $\square$

\noindent {\bf Corollary 2.2.} {\sl Let $M$ be a solid torus with a
triangulation $T$. Let $\gamma$ be a normal curve in $\partial M$
that bounds a meridian disc. Suppose that $\gamma$ has shortest length among
all such curves. Then $\gamma$ bounds a normal meridian disc.}

\noindent {\sl Proof.} Suppose that $\gamma$ intersects some edge
of $\partial M$ in two points of opposite sign. We may assume
that these two points of intersection are adjacent on the edge.
Let $\rho$ be the sub-arc of the edge between then. Then,
the interior of $\rho$ lies in the annulus $\partial M - \gamma$.
By our assumption about orientations, $\rho$ is inessential
in this annulus. Hence, a sub-arc of $\gamma$ is parallel in
$\partial M$ to $\rho$. We may isotope this sub-arc onto $\rho$,
and then reduce its length. This is contrary to hypothesis. 
So, the conditions of Proposition 2.1 are satisfied, and hence
$\gamma$ bounds a normal meridian disc. $\square$

We will now give the definition of a parallelity bundle. Let $M$ be a compact orientable
3-manifold with a triangulation $T$. Let $S$ be a properly embedded
normal surface in $M$, and let $X$ be the manifold obtained by cutting
$M$ along $S$. In other words, $X$ is the closure of $M - N(S)$.
Thicken the simplices of $T$ into handles, forming a handle structure $\hat {\cal H}$. 
Thus, each $i$-simplex of $T$ becomes an $i$-handle.  Then $S$ intersects
each handle of $\hat {\cal H}$ in a collection of discs. We say that two such
discs are of the same {\sl type} if there is an ambient isotopy of $M$,
preserving all the handles, that takes one disc to the other. If $D_1$ and
$D_2$ are adjacent discs of the same type in a handle $H$, then the component
of $H - {\rm int}(N(D_1 \cup D_2))$ that lies between them is
a product region $D^2 \times I$ where $D^2 \times \partial I$ is
parallel to $D_1 \cup D_2$. The union of these product regions is the
{\sl parallelity bundle} ${\cal B}$ in $X$. (See Figure 2.) The product structures patch together
to form an $I$-bundle structure on ${\cal B}$. 
Its {\sl horizontal boundary} $\partial_h {\cal B}$ is the $\partial I$-bundle,
and it lies in $\partial N(S)$.
The {\sl vertical boundary} $\partial_v {\cal B}$ is ${\rm cl}(\partial {\cal B}
- \partial_h {\cal B})$.
The vertical boundary $\partial_v {\cal B}$ of ${\cal B}$ is a collection of annuli. The
intersection of $\partial_v {\cal B}$ with $\partial M$ is a union of fibres in the $I$-bundle structure.
This is because this is true in each product region $D^2 \times I$ that makes up
${\cal B}$.

\vskip 18pt
\centerline{
\epsfxsize=3.5in
\epsfbox{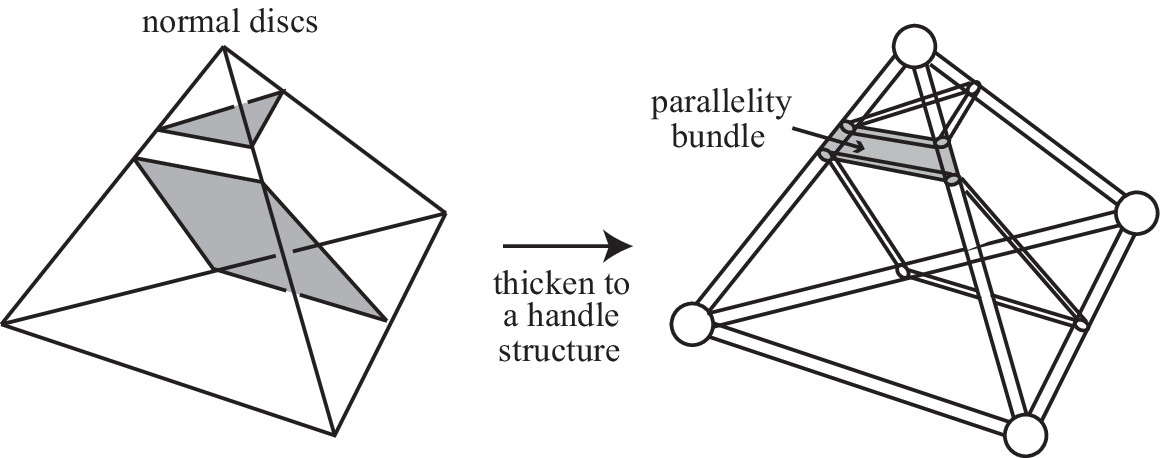}
}
\vskip 6pt
\centerline{Figure 2.}

\vskip 18pt
\centerline{\caps 3. Proof of the main theorems}
\vskip 6pt

In this section, we prove Theorems 1.1 and 1.2.

\noindent {\sl Proof.} 
Let $D$ be a properly embedded meridian disc, in normal form with respect to $T$.
Define the {\sl weight} of $D$ to be the number of points of intersection between
$D$ and the 1-skeleton of $T$. Define the {\sl length} of
$\partial D$ to be its number of points of intersection with the
1-skeleton. We define the {\sl complexity} of $D$ to be an ordered pair of non-negative
integers, the first of which is the length of $\partial D$, the
second of which is the weight of $D$. We order these
pairs lexicographically, which is a well-ordering, and choose $D$
so that it has minimal complexity. Note that this definition of complexity is a little
different from most treatments of the subject, where the weight of $D$ is given more
significance than the length of $\partial D$. (See [3] and [8] for example.) However, this alternative definition of
complexity is more suited to our purposes.

Note that $\partial D$ is a shortest
normal curve in $\partial M$ that bounds a meridian disc.
For suppose that there were a shorter normal meridian curve.
By Corollary 2.2, this would bound a normal meridian disc, which would therefore
have smaller complexity than $D$, which is a contradiction.

Let $N(D)$ be a thin regular neighbourhood of $D$ consisting of a union of normally parallel
copies of $D$. Let $X$ be the result of removing ${\rm int}(N(D))$ and ${\rm int}(N(D) \cap \partial M)$
from $M$. Thus, $X$ is a 3-ball, with two copies of $D$ in its boundary, which we denote by
$D_-$ and $D_+$. 
Let $A$ be the annulus $\partial X - {\rm int}(D_- \cup D_+) = {\rm cl}(\partial M - N(D))$.

Thicken each $i$-simplex of $T$ to an $i$-handle, forming a handle structure
$\hat {\cal H}$ for $M$. Let ${\cal B}$ be the parallelity bundle for $X$,
as described in Section 2. Let ${\cal B}'$ be the union of those components of ${\cal B}$ that
intersect $A$. Some possible configurations for components of ${\cal B}$ are shown in Figure 3.
Only the top two lie in ${\cal B}'$.

\vskip 18pt
\centerline{
\epsfxsize=4in
\epsfbox{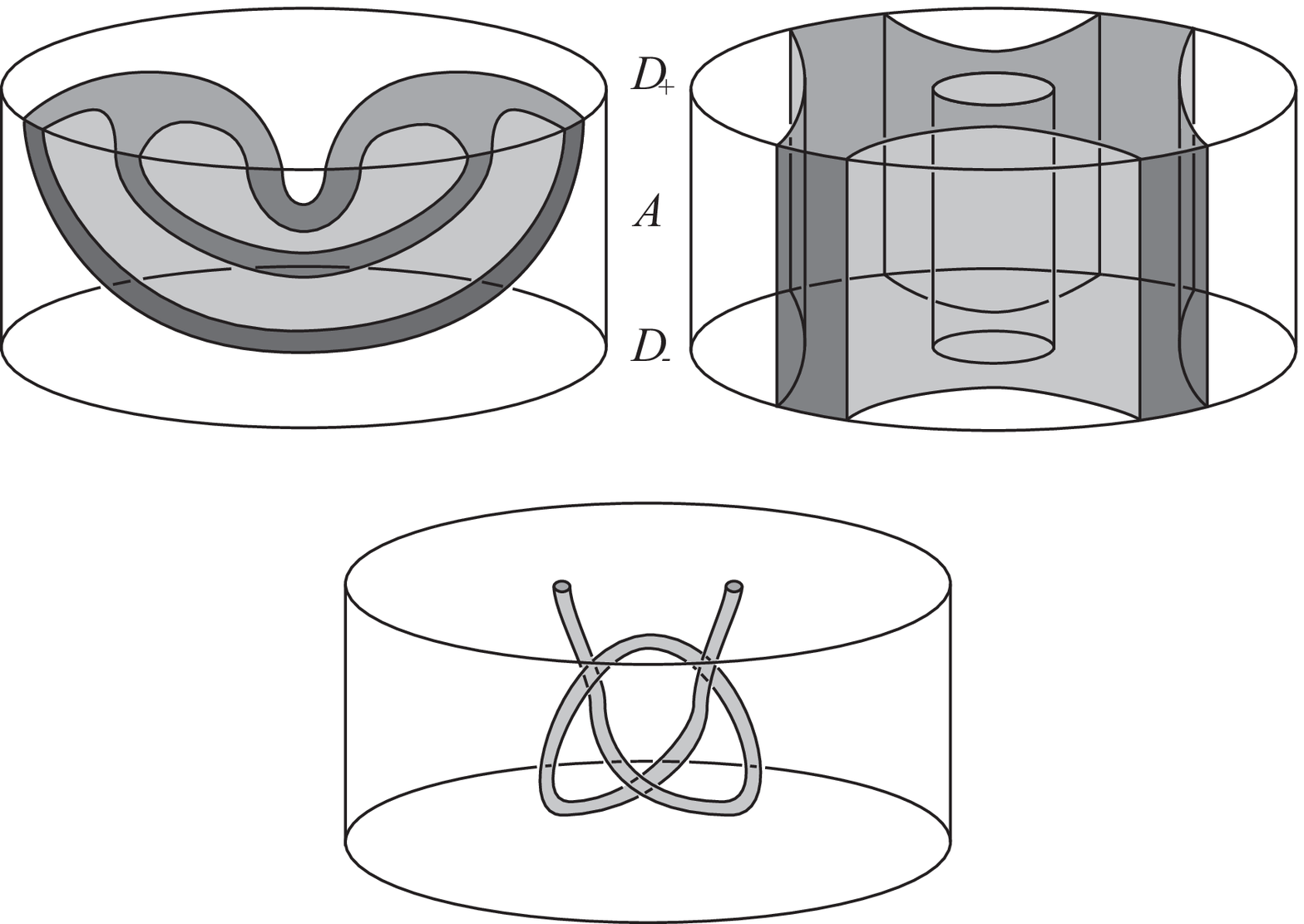}
}
\vskip 6pt
\centerline{Figure 3.}

\noindent {\sl Claim 1.} The $I$-bundle structure on ${\cal B}$ is in fact a
product structure.

Now, ${\cal B}$ is an $I$-bundle over a surface $F$, and to assert that ${\cal B}$
is a product is equivalent to asserting that $F$ is orientable, because the
total space of ${\cal B}$ is orientable. Suppose therefore that some component $F'$ of
$F$ is non-orientable. Let $E$ be the corresponding component of ${\cal B}$.
Then the horizontal boundary of $E$ is connected.
Double the 3-ball $X$ along $A$ to form a copy of $S^2 \times I$. Glue the
two copies of $E$ in $S^2 \times I$ together to form
a bundle $E^+$ in $S^2 \times I$. Now, $E^+$ need not be connected
(because $E$ may be disjoint from $A$), but if it is not,
discard one of its components. Because $E \cap A$ is a union of fibres, $E^+$ is an
$I$-bundle over a surface $F^+$. Note that $F^+$ is non-orientable,
since we are assuming that $F'$ is non-orientable. Hence, the horizontal boundary
of $E^+$ is connected. It is a planar surface lying in one component of
$S^2 \times \partial I$.
Each of its boundary components therefore bounds a disc with interior disjoint
from $E^+$. We can view $F^+$ as the zero section of the $I$-bundle
$E^+$. Each of its boundary components is a core curve of an annular component
of $\partial_v E^+$. It is therefore parallel to a boundary component
of this annulus, which in turn bounds a disc in $S^2 \times \partial I$.
Thus, we may extend $F^+$ to form a closed embedded non-orientable surface in $S^2 \times I$.
Then attaching 3-balls, we get such a surface in $S^3$. But it is well known that
$S^3$ contains no such surface. This proves the claim.

\noindent {\sl Claim 2.} Each component of ${\cal B}'$ 
must intersect both $D_-$ and $D_+$.

Suppose that, on the contrary, some component $E$ of ${\cal B}'$ 
is disjoint from $D_-$ or $D_+$, as in the top left of Figure 3, for example. Then each fibre
of $E \cap A$ is inessential in $A$. So some component of ${\rm cl}(A - E)$ is a disc $P$, with
boundary consisting of an arc $\alpha$ in $\partial A$ together with an $I$-fibre $\beta$ in
$\partial_v {\cal B}'$. Since $\beta$ is disjoint from 
the 1-skeleton of $T$, it therefore has
zero length. Now, $\alpha$ is parallel to a sub-arc $\alpha'$ in $\partial D$.
This has positive length, because the endpoints of $\beta$ lie in distinct normal discs of $D$.
We can attach $P$ to $D$, forming a new meridian disc for the solid torus $M$.
This has the effect of replacing $\alpha'$
with $\beta$, which reduces
the length of $\partial D$, and this is a contradiction.

\noindent {\sl Claim 3.} We may pick a product structure $D \times [-1,1]$ on $X$, so
that $D_- = D \times \{ -1 \}$, $D_+ = D \times \{ 1 \}$, and
so that the product structure agrees with that on ${\cal B}'$.

Let $S_+$ and $S_-$ be the surfaces $D_+ \cap {\cal B}'$
and $D_- \cap {\cal B}'$.
Pick a maximal collection $\alpha$ of disjoint non-parallel arcs properly embedded in
$S_+$, that are each essential in $S_+$
and that have boundary lying in $\partial D_+$. 
Because each component of ${\cal B}'$ intersects $\partial D_+$,
by the definition of ${\cal B}'$, each component of $S_+ - \alpha$ is therefore a
disc or an annulus. (Note that $\alpha$ may be empty, if $S_+$ is already 
a union of discs and annuli.) Let $\alpha'$ be the arc components of
${\rm cl}(\partial S_+ - \partial D_+)$. Using the
product structure on ${\cal B}'$, which exists by Claim 1, pick a collection of
embedded vertical discs $V$ in ${\cal B}'$ such that
$V \cap D_+ = \alpha \cup \alpha'$, and which are properly embedded
in $X$. Each such disc has boundary consisting of an arc in $D_+$,
an arc in $D_-$ and two arcs in $A$, by Claim 2. Cut $X$ along these discs, giving a collection of
balls. We will now pick a product structure on these balls,
which will patch together to give the required product
structure on $X$. Each ball $B$ intersects $D_+$ in a
single disc and intersects $D_-$ in a single disc.
Its intersection with $S_+$ is either empty, a disc or an
annulus. In the former case, we may clearly pick a
product structure on $B$ so that it agrees with product
structure on ${\cal B}' \cap B$. When $B \cap S_+$ is a
disc, $B$ lies in ${\cal B}'$, and so it already has
the desired product structure. When $B \cap S_+$ is an
annulus, one of its boundary components is disjoint
from $\partial D_+$. This bounds a disc in $D_+$.
Attached to this is a vertical annulus
in ${\cal B}'$. Its other boundary component lies in $D_-$,
and therefore bounds a disc in $D_-$. The union of the
discs in $D_-$ and $D_+$ with the vertical annulus bounds
a ball in $B$ with interior disjoint from ${\cal B}'$.
Hence, we may extend the product structure of ${\cal B}' \cap B$
over this ball, to give the required product structure on $B$.
This proves the claim.

Note that the claim is not obviously true with ${\cal B}'$ replaced
by ${\cal B}$. For example, there may be components of ${\cal B} - {\cal B}'$
which are disjoint from $D_-$ or $D_+$. Alternatively, they may lie in $X$ in a 
knotted way. (See the bottom of Figure 3.)

This product structure on $X$ determines a homeomorphism $\phi_0 \colon D_+ \rightarrow D_-$,
by first including $D_+$ into $D \times [-1,1]$ and then
projecting onto $D \times \{ -1 \} = D_-$. We also have a homeomorphism
$\psi \colon D_- \rightarrow D_+$ given by the gluing map.

Note that there was some flexibility in the choice of product structure on
$X - {\cal B}'$ and we will now vary it a little. In fact, it
is somewhat simpler to vary the homeomorphism $\phi_0$ to a new
homeomorphism $\phi \colon D_+ \rightarrow D_-$, which agrees with $\phi_0$ on $D_+ \cap {\cal B}'$
and which is isotopic to $\phi_0$ via an isotopy that is the
supported away from $D_+ \cap {\cal B}'$.

\vskip 12pt
\centerline{
\epsfxsize=3in
\epsfbox{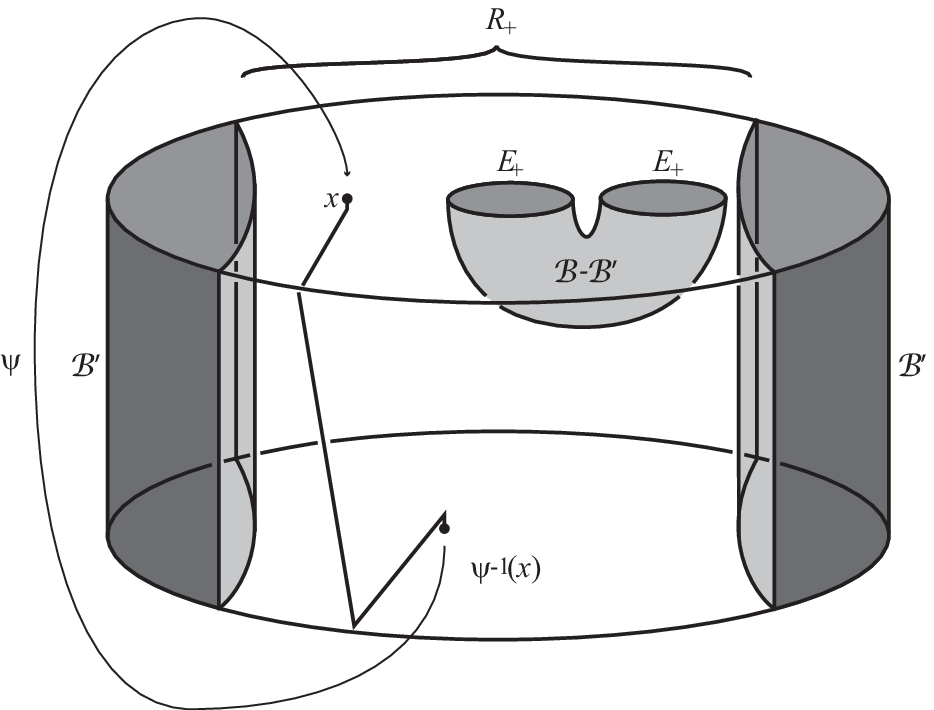}
}
\vskip 6pt
\centerline{Figure 4.}

\noindent {\sl Claim 4.} We may find such a homeomorphism $\phi$
with the following properties. Let $R_+$ be any component of 
$D_+ - {\cal B}'$ that intersects $\partial D_+$. Let $R_-$ be its image under $\phi$.
Let $E_+$ be the union of the discs in $R_+$ bounded by the curves
$R_+ \cap \partial_v({\cal B} - {\cal B}')$. (Possibly, $E_+$
is empty.) Similarly, let $E_-$ be the union of the discs in $R_-$ bounded by the curves
$R_- \cap \partial_v({\cal B} - {\cal B}')$. Then, we may ensure that:
\item{1.} $\phi(E_+)$ is disjoint from
$\psi^{-1}(E_+)$, and
\item{2.} $\phi^{-1}(E_-)$ is disjoint from
$\psi(E_-)$.

\noindent See Figure 4 for an example of $R_+$ and $E_+$.

We first ensure that (1) holds. Now,
$\psi^{-1}(E_+)$ is a collection of discs in the interior of $D_-$.
Hence, by isotoping the discs $\phi_0(E_+)$ sufficiently
close to $\partial D_-$, we may ensure that they are disjoint from
$\psi^{-1}(E_+)$. Thus, this gives  a homeomorphism $\phi_{00} \colon D_+ \rightarrow D_-$
such that $\phi_{00}(E_+)$ is disjoint from $\psi^{-1}(E_+)$.
(See Figure 5.)

\vskip 18pt
\centerline{
\epsfxsize=3.5in
\epsfbox{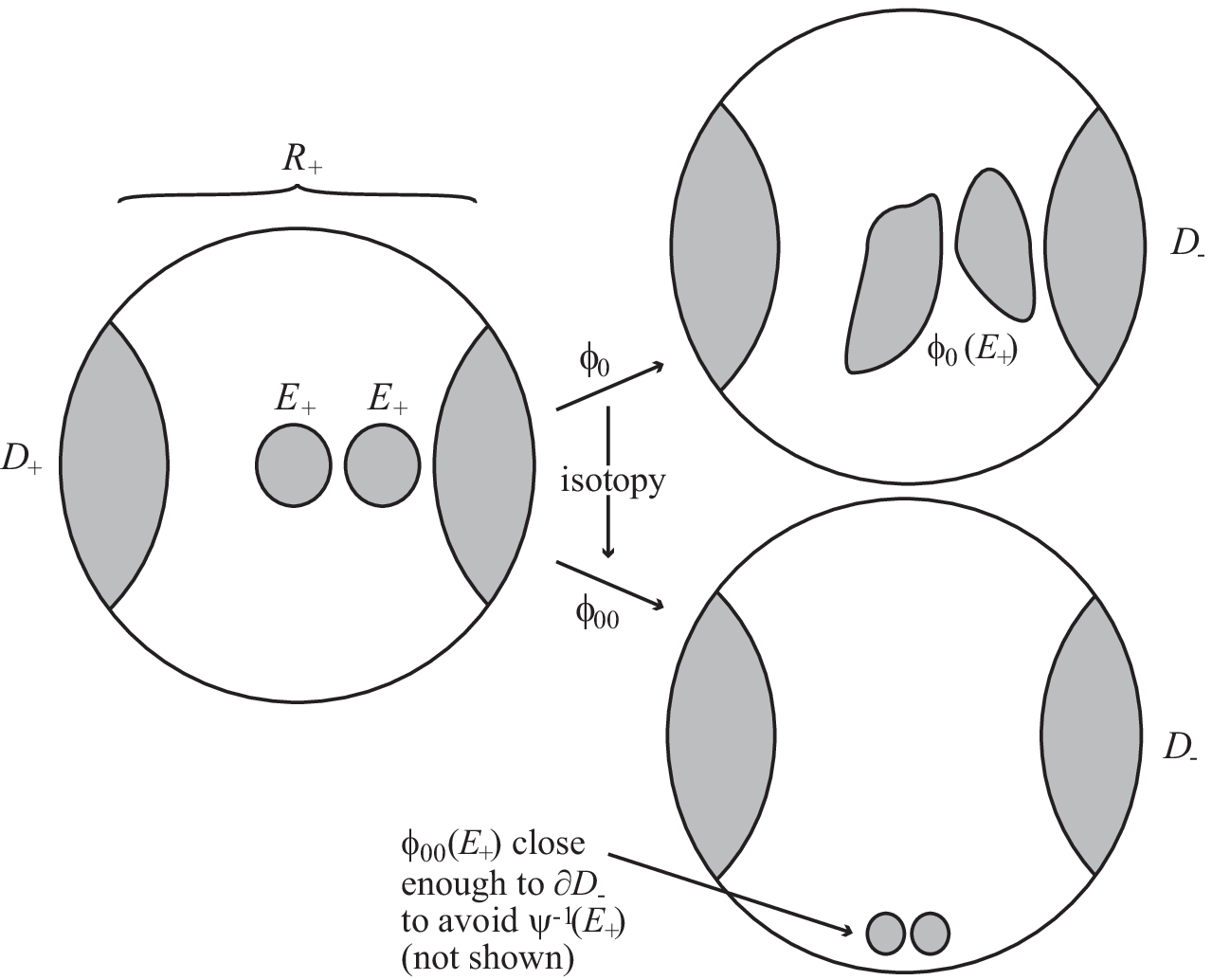}
}
\vskip 6pt
\centerline{Figure 5.}

We now ensure that (2) can also be achieved, without violating (1).
All we did was isotope $\phi_0$ to $\phi_{00}$ so that 
the discs $\phi_{00}(E_+)$ were close
to $\partial D_-$. Hence, we may also arrange that the discs $\phi_{00}(E_+)$ miss the
discs $E_-$. Therefore, $\phi_{00}^{-1}(E_-)$ is disjoint
from $E_+$. Now, there is a 1-parameter family of
homeomorphisms $D_+ \rightarrow D_+$, each supported in $R_+$,
starting at the identity, and which takes $\phi_{00}^{-1}(E_-)$
very close to $\partial D_+$. In particular, we may ensure that the image
of $\phi_{00}^{-1}(E_-)$ under the final homeomorphism is disjoint from
$\psi(E_-)$, because this is a collection of discs in the interior of $D_+$. Since $E_+$ does not
separate $\phi_{00}^{-1}(E_-)$ from $\partial D_+$, we may
assume that these homeomorphisms are the identity on
$E_+$. Thus, pre-composing these homeomorphisms with $\phi_{00}^{-1}$,
we obtain a 1-parameter family of homeomorphisms $D_- \rightarrow D_+$, starting at
$\phi_{00}^{-1}$, and which equal $\phi_{00}^{-1}$ throughout when restricted
to $\phi_{00}(E_+)$ and $D_- - R_-$. Inverting, we obtain a 1-parameter family of
homeomorphisms $D_+ \rightarrow D_-$, starting at $\phi_{00}$, and which equal $\phi_{00}$ throughout
when restricted to $E_+$ and $D_+ - R_+$. Let $\phi \colon D_+ \rightarrow D_-$ be the final
homeomorphism in this family. Since $\phi_{00}|E+ = \phi|E_+$,
condition (1) is therefore preserved. But we have ensured that
$\phi^{-1}(E_-)$ is disjoint from $\psi(E_-)$.
This gives condition (2), which proves the claim.

Thus, we have picked a homeomorphism $\phi \colon D_+ \rightarrow D_-$.
We also have a homeomorphism $\psi \colon D_- \rightarrow D_+$ arising from the
gluing map. Consider their composition $\psi \phi \colon D_+ \rightarrow D_+$.
By the Brouwer fixed-point theorem, $\psi \phi$ has a fixed point $x$ in $D_+$.
After an arbitrarily small perturbation of $\phi$, we may assume that
$x$ lies in the interior of a triangle or square of $D_+$.

\noindent {\sl Case 1.} $x \in {\cal B}'$.

Then the endpoints of the arc $x \times [-1,1]$ patch together to form
a core curve. But this leads to a contradiction in this case, because
the two endpoints of $x \times [-1,1]$ lie in distinct
normal discs of $D$.

\noindent {\sl Case 2.} $x$ lies in a component of $D_+ - {\cal B}'$ that
has non-empty intersection with $\partial D_+$.

Let $R_+$ be the component of $D_+ - {\cal B}'$ containing $x$, and let
$R_- = \phi(R_+)$. Define $E_+$ and $E_-$ as in Claim 4. An example of
such a configuration is shown in Figure 4.

Then $x$ is disjoint from ${\cal B}$ and lies in the component of $R_+ - {\cal B}$ that has non-empty
intersection with $\partial D_+$. For otherwise, $x$ lies in
$E_+$. By (1) in Claim 4, $\psi \phi(E_+)$ is disjoint from $E_+$. 
In particular, $\psi\phi(x) \not= x$, which is a contradiction.

Similarly, $\phi(x)$ is disjoint from ${\cal B}$ and lies in the component of $R_- - {\cal B}$
that has non-empty intersection with $\partial D_-$. For otherwise,
$\phi(x)$ lies in $E_-$. Hence, $x$ lies in $\phi^{-1}(E_-)$.
Also $\psi\phi(x)$ lies in $\psi(E_-)$. But, by (2) in Claim 4,
$\phi^{-1}(E_-)$ and $\psi(E_-)$ are disjoint.

We will shortly construct our pre-core curve $C$. We will pick an arc
in $A - {\cal B}$ that runs from $\partial R_+$ to $\partial R_-$.
We will then attach arcs in $D_+ - {\cal B}$ and $D_- - {\cal B}$ that run to $x$ and
$\phi(x) = \psi^{-1}(x)$ respectively. These three arcs patch together to
form an arc with one endpoint in $D_+$ and one in $D_-$. Gluing the
two ends of this arc will form the pre-core curve $C$. 
(See Figure 4.)

We first construct a pre-core curve $C_0$. The curves $C$ and $C'$ required by Theorems 1.1
and 1.2 will be minor modifications of this.
We wish to arrange that $C_0$ sits well with respect to $T$. 
Now, $D_+$, $D_-$ and $A$ clearly inherit cell structures from $T$. 
(For example, each 2-cell of $D_+$ and $D_-$ is a normal triangle or square.)
Note that whenever a $k$-cell of $A$ (respectively, $D_- \cup D_+$) lies in
${\cal B}$, then every incident $j$-cell of $A$ (respectively, $D_- \cup D_+$), with $j < k$,
also lies in ${\cal B}$. Thus, we may arrange that $C_0 \cap A$ 
intersects only the 2-cells and 1-cells and that it is transverse to the 1-cells. 
We ensure that $C_0 \cap A$ intersects
the 1-cells in as few points as possible, among all such arcs in $A - {\cal B}$ joining $R_+$
to $R_-$. Similarly, we may arrange that $C_0 \cap (D_- \cup D_+)$ 
intersects only the 2-cells and 1-cells of
these cell structures and that it is transverse to the 1-cells.
We choose these arcs $C_0 \cap (D_- - {\cal B})$ and $C_0 \cap (D_+ - {\cal B})$
to have shortest possible length among all arcs joining the endpoints of
$C_0 \cap A$ to $x$ and $\phi(x)$ satisfying the above conditions. 
This has various consequences:
\item{1.} The intersection of $C_0$
with each cell of $D_+$, $D_-$ and $A$ is connected whenever it is non-empty.
\item{2.} If $C_0$ intersects a 2-cell of $A$ that has non-empty
intersection with $R_+$ or $R_-$, then the relevant endpoint of $C_0 \cap A$
lies in that 2-cell. 
\item{3.} $C_0$ is disjoint from each 1-cell of $A$ that has an
endpoint in $R_- \cup R_+$ but does not lie entirely in $R_- \cup R_+$.

Now, $C_0$ is not in the form required by either Theorem 1.1 or 1.2, because
when it lies in $\partial M$, it lies in the 2-skeleton of $T$, but
when it runs through the interior of $M$, it is transverse to the 2-skeleton.
We therefore now explain how to modify $C_0$.

Our first aim is to create the core curve $C$. Let $C$ initially be $C_0$.
Now shorten $C \cap D_+$ a little at the $x$ end, so that
it terminates in the interior of a 1-cell of $D_+$. We correspondingly lengthen
$C \cap D_-$. We also modify $C \cap A$ near its endpoints so
that it ends on 0-cells of $A \cap (\partial D_- \cup \partial D_+)$. 
We now homotope $C \cap (D_- \cup D_+)$ into the 1-skeleton of $D_- \cup D_+$,
keeping the endpoints of these arcs fixed. 
This may create non-embedded arcs, but if so there is an obvious way to shorten
them. Thus, we may assume that $C \cap (D_- \cup D_+)$ is two embedded
arcs. Hence, $C$ is a pre-core curve. Note that $C \cap (D_- \cup D_+)$ need not be disjoint from
${\cal B}$. However, each 1-cell that it runs over is adjacent to a 2-cell
with interior disjoint from ${\cal B}$. 
Note that $C$ lies in the 2-skeleton of $T$, misses the vertices of $T$
and intersects the edges of $T$ in only finitely many points.
Thus, the following claim will prove Theorem 1.1.

\noindent {\sl Claim 5.} $C$ intersects the interior of each face of $T$ in at most
$10$ straight arcs.

The arcs of $C$ come in two types: parts lying in $A$, and parts lying in
$D_- \cup D_+$. By construction, both types of arc are straight. Thus, we only
need to bound the number of such arcs in any face of $T$. Consider any such
face $F$. The intersection $F \cap D$ consists of at most 3 types of normal
arcs. If a normal arc type does not arise, we add one in, for the sake
of streamlining the argument. By adding in such an arc very close to
a vertex of $F$, we can avoid it intersecting $C$. Thus, $F - N(D)$ consists of
three types of region:
\item{1.} triangular regions containing a vertex of $F$;
\item{2.} a hexagonal region containing all three arc types in its boundary;
\item{3.} rectangular regions between parallel arcs of $F \cap D$, which we call
{\sl parallelity rectangles}. 

The parallelity rectangles lie in ${\cal B}$, 
but $C$ need not be disjoint from ${\cal B}$.
However the parts of $C$ that lie in ${\cal B}$
are in $D_- \cup D_+$ and are adjacent to 2-cells of $D_- \cup D_+$
with interior disjoint from ${\cal B}$. Thus, if $C$ intersects a parallelity rectangle,
then in one of the adjacent 3-simplices, the parallelity rectangle lies
between a triangle and square of $D$. This can happen only once in
each of the adjacent 3-simplices. Thus, $C$ intersects at most 
$2$ parallelity rectangles, and when it intersects a parallelity
rectangle, it does so in at most $2$ arcs. It intersects each triangular region in
at most one arc. It intersects the hexagonal region in at most $3$ arcs.
So, the total number of arcs of $C \cap F$ is at most $10$,
which proves the claim.

We now explain how to construct the core curve $C'$ required
by Theorem 1.2. We simply push $C_0 \cap A$ a little into
the interior of $M$. This can be done in such a
way that $C'$ intersects each 3-simplex of $T$ in a collection of properly
embedded arcs.

\noindent {\sl Claim 6.} $C'$ intersects each 3-simplex of $T$ in at most $18$ arcs.

Let $\Delta$ be a 3-simplex of $T$, and let $P$ be the closure
of a component of $\Delta - N(D)$. Then, the boundary of $P$
consists of triangles and squares of $D_- \cup D_+$, together with bits
lying in $\partial \Delta$. The intersection $P \cap C_0$ is a union of
the following pieces:
\item{1.} arcs lying in $\partial \Delta$, which we call {\sl boundary} arcs;
\item{2.} arcs lying in a copy of a normal triangle or square of $D$, which
we call {\sl interior} arcs;
\item{3.} isolated points in edges of $\partial \Delta$ that lie in $\partial M$.

It is possible for interior arcs and boundary arcs to meet at their endpoints, 
and for boundary arcs to be joined together, but
otherwise these pieces are disjoint. (See Figure 6.)

Note that isolated points of $P \cap C_0$ arise as follows. When $C_0 \cap A$ runs over
an edge $e$ of $T$, all the adjacent 3-simplices of $T$ pick up a point of intersection
with $C_0$. If $\Delta$ is one of these 3-simplices and neither of the faces
of $\Delta$ containing $e$ lies in $\partial M$, this gives rise
to an isolated point of $\Delta \cap C_0$.

\vskip 18pt
\centerline{
\epsfxsize=2.5in
\epsfbox{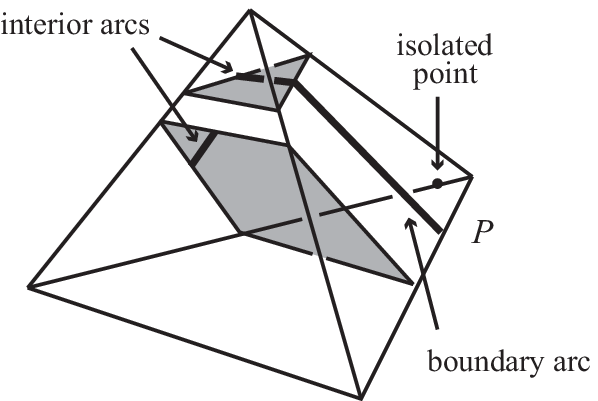}
}
\vskip 6pt
\centerline{Figure 6.}

Let $\Delta^1$ be the edges of $\Delta$, with the vertices removed.

In order to bound the number of components of $\Delta \cap C_0$,
the first thing that we do is remove all boundary arcs that lie in the same component
of $\Delta \cap C_0$ as an interior arc. This has the possible effect of increasing the number
of components of $\Delta \cap C_0$, but it does not decrease this number.
Once this has been done, every component of $\Delta \cap C_0$ is either an
interior arc or has non-empty intersection with $\Delta^1$.
Thus, the number of components of $\Delta \cap C_0$ is at most
the number of interior arcs plus the number of intersections between
$C_0$ and $\Delta^1$.

Now, each edge of $\Delta$ is divided into arcs by $D$.
If one of these arcs is disjoint from the vertices of $\Delta$,
then it lies in ${\cal B}$. Thus, the only components of $\Delta^1 - N(D)$ that
can intersect $C_0$ are those that are incident to a vertex of $\Delta$.
Call such components of $\Delta^1 - N(D)$ {\sl vertex-incident}.
There are at most $12$ vertex-incident arcs of $\Delta^1 - N(D)$,
and so there are at most 12 points of $C_0 \cap \Delta^1$.

Between adjacent normally parallel triangles or squares
of $\Delta \cap D$ lies ${\cal B}$. Hence, $C_0$ misses this region. So, each
normal disc type of $\Delta \cap D$ can give rise to at most two interior arcs.
There are at most $5$ normal disc types in $\Delta$, and so there
at most $10$ interior arcs of $C_0 \cap \Delta$.

Thus, the total number of components of $C_0 \cap \Delta$ is at most
$12 + 10 = 22$. But we can reduce this bound down to $18$ as follows.

Each interior arc of $\Delta \cap C_0$ lies in a normal disc $E$
of $D_- \cup D_+$. Then $E$ lies in $R_- \cup R_+$,
and so the components of $\Delta^1 - N(D)$ incident to $E$ are disjoint
from $C_0$, by condition (3) above. When such a component of
$\Delta^1 - N(D)$ is also vertex-incident, this
reduces the number of possible points of $C_0 \cap \Delta^1$
below $12$. Now, $6$ discs of $(D_- \cup D_+) \cap \Delta$ that are
disjoint from ${\cal B}$ may be disjoint from vertex-incident
arcs of $\Delta^1 - N(D)$. But if there are any more than
$6$ interior arcs of $\Delta \cap C_0$, then each one reduces
the number of possible points of $\Delta^1 \cap C_0$ by at least $3$.
So, we deduce that the total number of components of
$C_0 \cap \Delta$ is at most $12 + 6 = 18$. This proves the claim.

\noindent {\sl Claim 7.} $C' \cap \Delta$ is parallel to 
a collection of arcs $\alpha$ in $\partial \Delta$, with
the property that the intersection between each component of $\alpha$
and each face of $\Delta$ is at most one straight arc.

We have already arranged that $C'$ intersects each 3-simplex in one of only
finitely many possible configurations. So, one could simply perform a case-by-case
check to prove this claim. However, there are many cases, and so we present a
uniform argument.

Let $\Delta$ be a 3-simplex of $T$, and let $P$ be the closure of a component
of $\Delta - N(D)$. We first show that $C' \cap P$ is parallel in $P$ to a collection
of arcs $\alpha$ in $\partial \Delta \cap P$. We will then show that $\alpha$ can
be chosen so that each component intersects each face of $\Delta$ in at most one straight arc.
This will prove the claim.

Now, the boundary arcs and isolated points of $C_0 \cap P$ already lie
in $\partial \Delta$. These will form part of $\alpha$. Thus, we need only
to worry about the parts of $C_0 \cap P$ lying in interior arcs. There is at
most one interior arc lying in each normal triangle or square of $D_- \cup D_+$, and this can be slid,
keeping its endpoints fixed, into the boundary of the triangle or square.
There are two possible directions in which this slide can be performed,
and we choose a direction which minimises the number of intersections with
$\Delta^1$. This creates our choice of arcs $\alpha$.

How could a component of $\alpha$ intersect more than one face of $\Delta$?
Each component of $C_0 \cap \Delta$ intersects each face in at most one
component. Thus, the only possible problem that could arise is when
the process of sliding creates a new component of intersection with
a face. This only happens when the interior arc of $C_0 \cap \Delta$
lies in a normal square $E$ and the interior arc joins opposite edges
of the square. Then one of the sides of the square ($\beta$, say)
becomes a new component of intersection between $\alpha$ and a face
of $\Delta$. Now, each component of $\Delta - D$ that does not lie wholly in
${\cal B}$ contains at most one normal square of $D_- \cup D_+$. So, the only
way that a problem could arise is when the 2-cell of $A$
containing $\beta$ has non-empty intersection with $C_0$. 
But in this case, condition (2) above implies that the relevant
endpoint of $C_0 \cap A$ lies in that 2-cell. This can happen in
only that 2-cell. So, if we slide $C_0 \cap E$ the other
way across $E$, then this problem is avoided. Finally,
note that $\alpha$ intersects each face of $P$ in at most one arc.
We may therefore isotope $\alpha$ so that it is straight in each face.
This proves the claim.

Thus, we have proved Theorems 1.1 and 1.2 in this case.

\noindent {\sl Case 3.} $x$ lies in a component of $D_+ - {\cal B}'$
that is disjoint from $\partial D_+$.

An example of such a configuration is shown in Figure 7.

\vskip 18pt
\centerline{
\epsfxsize=2.5in
\epsfbox{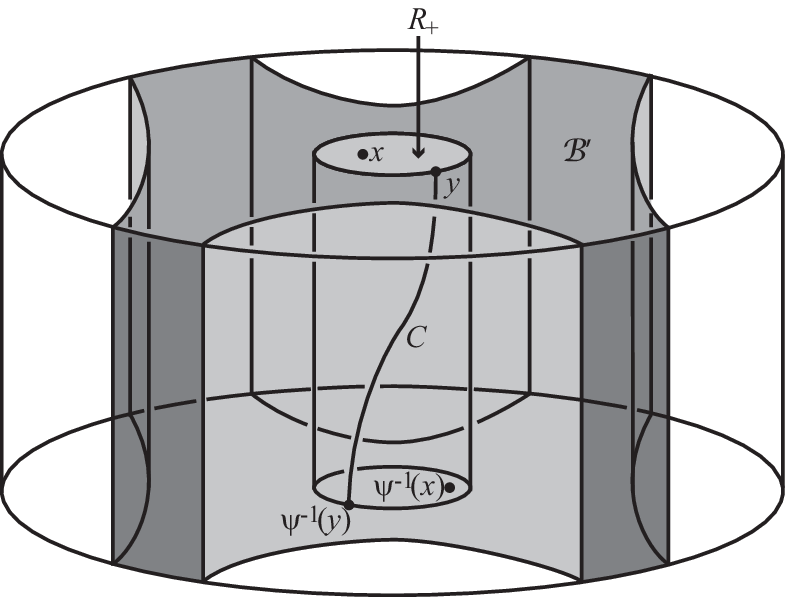}
}
\vskip 6pt
\centerline{Figure 7.}

Let $R_+$ be the closure of the component of $D_+ - {\cal B}'$ containing $x$.
This is a disc disjoint from $\partial D_+$.
Since $\psi \phi (x) = x$, $\psi \phi(R_+)$ has non-empty intersection
with $R_+$. We consider two cases.

\noindent {\sl Case 3A.} $\psi \phi(\partial R_+) \cap \partial R_+ \not= \emptyset$.

Thus, $\psi \phi (\partial R_+) \cap \partial R_+$ contains some point $y$.
If we were to pick an arc in the annulus $\partial R_+ \times [-1,1]$ running from $y$
to $\psi^{-1}(y)$, then this would patch up to form a core curve for $M$.
In fact, we will feel free to vary this arc by an ambient isotopy, keeping its
endpoints fixed. The arc first runs from $y$ to $\phi(y)$ vertically
across $\partial_v {\cal B}$. It then runs from $\phi(y)$ to $\psi^{-1}(y)$
along 1-cells of $D_-$. Each such 1-cell is adjacent to a 2-cell with
interior disjoint from ${\cal B}$. So, as in Case 2,
the resulting core curve $C$ lies in the 2-skeleton of $T$
and intersects the interior of each face of $T$ in at most $10$ straight arcs.
One may also perturb $C$ to create a core curve $C'$ that is
transverse to the 2-skeleton of $T$. It then intersects each tetrahedron
only in interior arcs, as described in Case 2. So, $C'$ intersects
each tetrahedron in at most $10$ arcs, and these satisfy the
conclusions of Theorem 1.2. Thus, the theorems are proved in this
case.

\noindent {\sl Case 3B.} $\psi \phi(\partial R_+) \cap \partial R_+ = \emptyset$.

Now, $\psi \phi(R_+)$
has non-empty intersection with $R_+$, and yet their boundaries
are disjoint, and so they must be nested. Say that $\psi \phi(R_+)$
lies in $R_+$. 

Note that $D_+$ inherits a handle structure from $\hat {\cal H}$.
Both $R_+$ and $\psi \phi(R_+)$ are a union of handles in this
handle structure. However, $R_+$ and $\psi \phi(R_+)$ do not
inherit handle structures. This is because it need not be true
that whenever an $i$-handle lies in $R_+$ (or $\psi \phi(R_+)$),
then so does every $j$-handle to which it is attached, with $j < i$.
In fact, ${\rm cl}(D_+ - R_+)$ and ${\rm cl}(D_+ - \psi \phi(R_+)$
inherit handle structures. This is because the handles of $D_+$
that intersect $R_+$ but do not lie in it are actually part of
${\cal B}$. And whenever an $i$-handle lies in ${\cal B}$,
then so is every $j$-handle to which it is attached, with $j < i$.

We now consider the weight of $R_+$ and $\psi \phi (R_+)$. Recall
that this is just the number of intersections with the $1$-skeleton
of the triangulation. Equivalently, it is the number of $0$-handles
of $R_+$ (or $\psi \phi(R_+)$).

\noindent {\sl Claim 8.} $\psi \phi(R_+)$ has strictly smaller weight than $R_+$.

Now, $\psi \phi(R_+)$ is a subset of $R_+$ consisting of a union of handles.
Thus, if the weight of $\psi \phi(R_+)$ is not less than that of $R_+$,
then only 1-handles and 2-handles are removed when constructing $\psi \phi(R_+)$
from $R_+$. But $\psi\phi(\partial R_+)$ is disjoint from $\partial R_+$,
and so every handle of $R_+$ that contains $\partial R_+$ is removed.
Also, if an $i$-handle is removed, then so is every $j$-handle that it is
attached to, whenever $j < i$. Now, if there is no 0-handle of $R_+$
that is incident to a 1-handle or 2-handle which contains an arc of $\partial R_+$,
then $R_+$ consists only of 1-handles and 2-handles. But in this case,
every handle of $R_+$ has non-empty intersection with $\partial R_+$,
and so all of $R_+$ is removed when creating $\psi \phi(R_+)$. In other
words, $\psi \phi(R_+)$ is empty, which is a contradiction.
This proves the claim.

So, consider the disc $D'$ obtained from $D$ by removing the interior of
$R_+$, and then attaching $(\partial R_+ \times [-1,1]) \cup \psi \phi(R_+)$.
Since $\partial R_+ \times [-1,1]$ is vertical in the parallelity
bundle, it has zero weight. Thus, $D'$ has the
same boundary as $D$ but has smaller weight. This contradicts our
assumption that $D$ has minimal complexity. $\square$

\vfill\eject
\centerline{\caps 4. Partially ideal triangulations and affine handle structures}
\vskip 6pt

In this section, we generalise Theorems 1.1 and 1.2 from triangulations to two other
types of representation of a 3-manifold: partially ideal triangulations and
affine handle structures.

A {\sl partially ideal triangulation} of a 3-manifold $M$ is an
expression of $M - \partial M$ as a collection of 3-simplices with their faces
identified in pairs, and then with some of their vertices removed.
We will assume, without loss of generality, that all of the gluing maps
are affine. Hence, as in previous sections, it makes sense to speak of
a straight line in a face of the triangulation.

One can remove a small open product neighbourhood of each end of $M - \partial M$,
thereby truncating the ideal vertices of each tetrahedron. This can be made into a
convex polyhedron is by starting with a regular Euclidean tetrahedron and truncating
some its vertices. We call this a {\sl truncated partially ideal triangulation}.

We will prove the following version of Theorem 1.2 in this context.

\noindent {\bf Theorem 4.1.} {\sl Let $T$ be a truncated partially ideal triangulation
of the solid torus $M$. Then there is a core curve $C$ of $M$ that intersects 
each truncated tetrahedron $\Delta$ of $T$ in a collection of at most $48$ properly embedded arcs, with
endpoints in the interiors of the faces of $\Delta$. These
arcs $C \cap \Delta$ are simultaneously parallel to a union of arcs in $\partial \Delta$,
each of which intersects each face in at most $6$ straight lines.}

We will also deal with handle structures.
Whenever we refer to a handle structure on a 3-manifold, we insist that each
handle is attached to handles of strictly lower index.
An {\sl affine handle structure} on a 3-manifold $M$ is a handle structure where each 
0-handle and 1-handle is identified with a compact (but possibly non-convex) polyhedron in ${\Bbb R}^3$, 
so that 
\item{1.} each face of each polyhedron is convex;
\item{2.} whenever a 0-handle and 1-handle intersect, each component of intersection is
identified with a convex polygon in ${\Bbb R}^2$, in such a way that the inclusion
of this intersection into each handle is an affine map onto a face of the
relevant polyhedron;
\item{3.} for each 0-handle $H_0$, each component of intersection with a 2-handle,
3-handle or $\partial M$ is a union of faces of the polyhedron associated with $H_0$;
\item{4.} the polyhedral structure on each 1-handle is the product of a convex 2-dimensional
polygon and an interval.

Here, we have the following version of Theorem 1.2.

\vfill\eject
\noindent {\bf Theorem 4.2.} {\sl Let ${\cal H}$ be an affine handle structure of the
solid torus $M$. Suppose that each $0$-handle of ${\cal H}$ has at most $4$ components
of intersection with the $1$-handles, and that each $1$-handle has at most $3$
components of intersection with the $2$-handles. Then $M$ has a core curve that intersects 
only the $0$-handles and $1$-handles, that respects the product structure on the
$1$-handles, that intersects each $1$-handle in at most $24$ straight arcs,
and that intersects each $0$-handle in at most $48$ arcs.
Moreover, the
arcs in each $0$-handle are simultaneously parallel to a collection of arcs $\alpha$ in the boundary of the
corresponding polyhedron, and each component of $\alpha$ intersects each face of the polyhedron
in at most $6$ straight arcs.}

The definition of an affine handle structure may not seem to be particularly natural.
For example, it is not clear why one only identifies the 0-handles and 1-handles with
polyhedra. However, to impose a polyhedral structure on the handles with
higher index would be unduly restrictive and also unnecesssary. Similarly, it is
necessary for us to insist that the faces of each polyhedron are convex, but
not that the polyhedra themselves are convex. 

In Theorem 4.2, the restrictions on the number of components of intersection between various types of
handle may also seem strange. But the conditions are natural. For example, they are satisfied by the handle structure
that is dual to a partially ideal triangulation. In fact, Theorem 4.2 is a generalisation of
Theorem 4.1. For, suppose that we are given a truncated partially ideal triangulation of a 3-manifold. 
We then declare that each truncated tetrahedron is a 0-handle. 
Each face of the ideal triangulation is made
into a 1-handle. Each edge becomes a 2-handle, and each non-ideal vertex becomes a
3-handle. Theorem 4.2 then provides a core curve that intersects this handle
structure nicely. Translating this back to the initial ideal triangulation,
we obtain the conclusion of Theorem 4.1.

We therefore focus on the proof of Theorem 4.2. Let ${\cal H}$ be the given
handle structure, and for $i = 0, 1, 2, 3$, let ${\cal H}^i$ be the
union of the $i$-handles. 

We would like to be able to assume that
${\cal H}^0 \cap \partial M$ is a collection of discs. However, this
need not be the case. But, suppose that, for some 0-handle $H_0$, 
$H_0 \cap \partial M$ has a non-disc component. Pick a simple closed curve
that is essential in $H_0 \cap \partial M$. We may find
such a curve that intersects each convex face of $H_0 \cap \partial M$
in at most one straight arc. This curve bounds
a properly embedded disc $E$ in $H_0$. There are two cases to
consider: where $E$ is boundary parallel in $M$, and where
$E$ is a meridian disc for $M$. If $E$ is boundary parallel, then
the parallelity region is a 3-ball. Remove from ${\cal H}$ all the
handles that lie in this ball, except $H_0$. It is clear that if
the theorem holds for this new handle structure, then it also holds
for ${\cal H}$. Thus, we may assume that this case does not arise.
Suppose now that $E$ is a meridian disc. Then the result of cutting
$\partial M$ along $\partial E$ is an annulus. If we pick any properly
embedded arc in this annulus joining its two boundary components, 
then this can be closed up to form a pre-core curve $C$. We may clearly
arrange that $C$ intersects only the 0-handles and 1-handles
of ${\cal H}$, and that it respects the product structure on the
1-handles. We may also ensure that it runs over each component of
${\cal H}^1 \cap \partial M$ at most once. Thus, it runs over each
1-handle of ${\cal H}$ at most three times. The points where it 
enters and leaves any 0-handle lie in ${\cal H}^1$, and
so it intersects each 0-handle in at most $6$ arcs. We may
ensure that each such arc intersects each face of $\partial {\cal H}^0$
in at most one straight arc. Now push $C$ a little into the
interior of $M$ to create the core curve required by Theorem 4.2.

We may therefore assume that ${\cal H}^0 \cap \partial M$ is
a collection of discs. 
This forces each 2-handle to run over at least one 1-handle.
Therefore, ${\cal H}^0 \cap ({\cal H}^1 \cup {\cal H}^2)$ forms a thickened
graph, where ${\cal H}^0 \cap {\cal H}^1$ is the thickened vertices
and ${\cal H}^0 \cap {\cal H}^2$ is the thickened edges. 

In the proof of Theorems 1.1 and 1.2, a handle structure was constructed from the triangulation of $M$,
and this was used to define the parallelity bundle. We want to do something similar here. The result will be a handle structure $\hat {\cal H}$
which is, in some sense, dual to the given handle structure ${\cal H}$. Handles of $\hat {\cal H}$
will arise in two possible ways. Each $i$-handle of ${\cal H}$ will give rise to a $(3-i)$-handle
of $\hat {\cal H}$. But there is a second type of handle of $\hat {\cal H}$. The handle structure
${\cal H}$ induces a handle structure for $\partial M$, and each $i$-handle of $\partial M$
gives rise to $(2-i)$-handle of $\hat {\cal H}$. One may also form a cell structure for $M$
by collapsing each $i$-handle of $\hat {\cal H}$ to an $i$-cell.

The first step in the proof of Theorems 1.1 and 1.2 was to find a meridian
disc $D$ in normal form. The complexity of $D$ was defined to be an ordered pair
of integers: the length of $\partial D$, and the weight of $D$. Complexities
were ordered lexicographically, and $D$ was chosen to have smallest possible
complexity. We perform a similar process in this case.

We start by making $D$ {\sl standard} in the handle structure ${\cal H}$. This means that
it intersects each handle in a collection of properly embedded discs, misses the 3-handles,
and respects the product structure on each 2-handle and 1-handle. We define the {\sl boundary weight}
of $D$ to be the number of intersections between
$D$ and $\partial M \cap {\cal H}^0$. We define
the {\sl interior weight} to be the number of intersections
between $D$ and ${\cal H}^2$. We define
the compexity of $D$ to be the ordered pair, boundary weight then interior weight,
and choose $D$ so that it has minimal complexity.

The resulting disc $D$ need not be normal, in the usual sense of the word
(for example, as in [3] or [5]). This is because there may be arcs of
intersection between $D$ and a component of ${\cal H}^0 \cap {\cal H}^1$
which run from a component of ${\cal H}^0 \cap \partial M$ to
an adjacent component of ${\cal H}^0 \cap {\cal H}^2$. But $D$ satisfies
the other conditions of normality:
\item{1.} $D \cap {\cal H}^0 \cap {\cal H}^1$ consists of arcs;
\item{2.} $D \cap {\cal H}^0 \cap {\cal H}^2$ consists of arcs respecting
the product structure in ${\cal H}^0 \cap {\cal H}^2$;
\item{3.} $D \cap {\cal H}^0 \cap \partial M$ consists of arcs;
\item{4.} no arc of intersection between $D$ and ${\cal H}^0 \cap {\cal H}^1$
has endpoints in the same component of ${\cal H}^0 \cap {\cal H}^1 \cap {\cal H}^2$,
or in the same component of ${\cal H}^0 \cap {\cal H}^1 \cap \partial M$;
\item{5.} no arc of intersection between $D$ and ${\cal H}^0 \cap \partial M$
has endpoints in the same component of ${\cal H}^0 \cap {\cal H}^1 \cap \partial M$;
\item{6.} each component of $D \cap \partial {\cal H}^0$ intersects each
component of ${\cal H}^0 \cap {\cal H}^2$ in at most one arc;
\item{7.} each component of $D \cap \partial {\cal H}^0$ intersects each
component of ${\cal H}^0 \cap \partial M$ in at most one arc.

For any handle $H$ of ${\cal H}$,
we say that two discs of $H \cap D$ are of the same {\sl type} if there is an
ambient isotopy, preserving all the handles of ${\cal H}$, that takes
one disc to the other.
The above conditions imply that $D$ intersects each handle in only finitely many disc types. 
By condition (6), each disc $E$ of $D \cap {\cal H}^0$ intersects each component of
${\cal H}^0 \cap {\cal H}^2$ in at most $1$ arc. By condition (7),
$E$ intersects each component of 
${\cal H}^0 \cap \partial M$ in at most $1$ arc.
By a combination of conditions (4), (6) and (7) and the assumption that 
each $1$-handle has at most $3$ components of intersection with the $2$-handles,
$E$ intersects each component of
${\cal H}^0 \cap {\cal H}^1$ in at most $3$ arcs.

We now make $D$ sit nicely with respect to $\hat {\cal H}$. Each component of
intersection between $D$ and a handle of ${\cal H}$ or $\partial M$
gives rise to a disc component of intersection between $D$ and a handle
of $\hat {\cal H}$. We define the ${\sl type}$ of such a disc
just as in the case of ${\cal H}$. Between adjacent discs of
the same type, there is a product region, and these patch together
to form the {\sl parallelity bundle} ${\cal B}$.

Now, each 0-handle $H_0$ of ${\cal H}$ has been identified with a polyhedron,
and each component of $H_0 \cap {\cal H}^1$, $H_0 \cap {\cal H}^2$
and $H_0 \cap \partial M$ is a union of convex faces. We may clearly
arrange that $D$ intersects each such convex face in a collection of
straight arcs. Moreover, we may ensure that each component of $D \cap H_0 \cap {\cal H}^1$,
$D \cap H_0 \cap {\cal H}^2$ and $D \cap H_0 \cap \partial M$ intersects each
convex face in at most one such arc.

The rest of the argument proceeds as in the case of triangulations,
apart from the exact specification of the core curve $C$ with respect to ${\cal H}$. 
The symbols $D_-$, $D_+$, $A$, $X$ and ${\cal B}'$ denote exactly
what they did before. A homeomorphism $\phi \colon D_+ \rightarrow D_-$
is picked, and $\psi \colon D_- \rightarrow D_+$ is the gluing map.
The composition $\psi \phi \colon D_+ \rightarrow D_+$ has a fixed point $x$.
The proof then divides into cases according to the location of $x$.
As before, there are two cases when a copy of $C$ is
constructed (Cases 2 and 3A). We will focus on Case 2, because this is
more complex. Recall that this is when $x$ lies in a component 
of $D_+ - {\cal B}'$ that has non-empty intersection with $\partial D_+$.
Let $R_+$ be the closure of the component of $D_+ - {\cal B}'$ containing $x$, and let
$R_- = \phi(R_+)$. Recall that in this case, we picked a curve $C_0$ in $M$,
avoiding ${\cal B}$, as follows. The arc $C_0 \cap A$ ran between
$R_-$ and $R_+$. It was disjoint from the 0-cells of $A$ and was transverse to the
1-cells. It was chosen to have fewest number of intersections with the
1-cells. Then $C_0 \cap (D_+ - {\cal B})$ and $C_0 \cap (D_- - {\cal B})$ were arcs running from the
endpoints of $C_0 \cap A$ to $x$ and $\phi(x)$. These were also chosen
to have smallest length among all curves in $D_{\pm} - {\cal B}$ joining these
specified endpoints,  in the sense of having fewest number of points of
intersection with the 1-cells.

In our situation, we do the same. Once again, $A$, $D_-$ and $D_+$ inherit cell
structures. For example, each component of intersection between $D$ and
${\cal H}^0$ gives a component of intersection between $D$ and a 3-handle of
$\hat {\cal H}$, and this becomes a 2-cell of $D_-$ and $D_+$. 
The arc $C_0 \cap A$ is chosen to
have smallest length among all curves in $A - {\cal B}$ joining $R_-$ to $R_+$.
We then pick shortest arcs in $(D_- - {\cal B})$ and $(D_+ - {\cal B})$ joining the endpoints of
$C_0 \cap A$ to $x$ and $\phi(x)$. In the proof of Theorem 1.2, we needed
to push $C_0$ a little into the interior of $M$ to form a core curve.
We do the same here, and let $C$ be the resulting core curve.

Let $H_0$ be a 0-handle of ${\cal H}$. We wish to find an upper bound on the
number of arcs of $C_0 \cap H_0$. To do this, we note that the only places where
$C_0$ can enter $H_0$ are in the discs $H_0 \cap {\cal H}^1$. So, we will bound
the number of points of $C_0 \cap H_0 \cap {\cal H}^1$. These points come in
two types: those lie that lie on the boundary of the discs $H_0 \cap {\cal H}^1$,
and those that lie in the interior of these discs. The points that lie
on the boundary of these discs lie in $\partial M$. So consider an arc
component of $H_0 \cap {\cal H}^1 \cap \partial M$. This is divided up by
$D$, and between adjacent points of $D$, there lies ${\cal B}$, which $C_0$
avoids. So, $C_0$ can intersect each arc of $H_0 \cap {\cal H}^1 \cap \partial M$
in at most $2$ points. There are at most $3$ such arcs in each component of
$H_0 \cap {\cal H}^1$, and so this
gives at most $6$ endpoints of $C_0 \cap H_0$. The second type of endpoint
of $C_0 \cap H_0$ lies in the interior of the discs of $H_0 \cap {\cal H}^1$.
These lie on an arc of $D \cap H_0 \cap {\cal H}^1$. Between two such
arcs of the same type, there again lies the parallelity bundle ${\cal B}$.
So, each arc type of $D \cap H_0 \cap {\cal H}^1$ gives at most two endpoints
of $C_0 \cap H_0$. There are at most $9$ arc types in each component of
$H_0 \cap {\cal H}^1$. So, this gives at most $18$ endpoints of $C_0 \cap H_0$. So, in total,
we have at most $6 + 18 = 24$ endpoints of $C_0 \cap H_0$ in each
component of $H_0 \cap {\cal H}^1$, and hence
at most $48$ arcs of $C_0 \cap H_0$.

We now need to justify why $C_0 \cap {\cal H}^0$ is parallel to a collection
of arcs in $\partial {\cal H}^0$ as described in Theorem 4.2.
Some of $C_0$ already lies in $\partial {\cal H}^0$, but we slide the
parts lying in $D \cap {\cal H}^0$ into $\partial (D \cap {\cal H}^0)$. 
The result is a collection of arcs $\alpha$ in $\partial {\cal H}^0$
to which $C_0 \cap {\cal H}^0$ is parallel.

We must show that $\alpha$ can be isotoped, keeping $\partial \alpha$ fixed,
so that it intersects each face of each polyhedron of ${\cal H}^0$
in a collection of straight arcs, and that each component of
$\alpha$ intersects each such face in at most $6$ arcs.
Now, $\partial {\cal H}^0$ is a
union of convex faces, and these have been divided up by $\partial D$ into
convex polygons. Thus, we may isotope $\alpha$ so that it intersects
each face in straight arcs. Now, $C_0$
is a concatenation of {\sl interior arcs}, which run over a component of $D \cap {\cal H}^0$
and {\sl boundary arcs}, which run over a component of $\partial M \cap {\cal H}^0$.
An interior arc is joined to a boundary arc exactly where $C_0 \cap A$ is joined
to $C_0 \cap (D_- \cup D_+)$, and there are just two of these points.
So, in each component of $C_0 \cap {\cal H}^0$, at most one boundary arc is used,
and at most two interior arcs. Now, each interior arc has been slid into the boundary
of a component of $D \cap {\cal H}^0$. Each component of $D \cap {\cal H}^0$
intersects each face of $\partial {\cal H}^0$ at most $3$ times. In fact,
if the face lies in $\partial M$, then the component of $D \cap {\cal H}^0$
runs over it at most once. Each boundary arc lies in a single component
of $\partial M \cap {\cal H}^0$, and so runs over each face at most once.
Thus, each component of $\alpha$ runs over each face of
$\partial {\cal H}^0$ at most $6$ times, as required.
$\square$

\vskip 18pt
\centerline{\caps 5. Riemannian metrics on solid tori}
\vskip 6pt

In this section, we will prove that any Riemannian metric on a solid torus $M$
with bounded sectional curvature and a lower bound on injectivity radius has a core curve
with length that is linearly bounded by the volume of $M$.

\noindent {\bf Theorem 1.3.} {\sl For each $K$, $I > 0$, there is a constant
$c(K,I)$ with the following property. If $M$ is a solid torus with a Riemannian
metric having volume at most $V$, injectivity radius at least $I$ and all sectional
curvatures in the interval $(-K,K)$, then $M$ contains a core curve with length
at most $c(K,I) \ V$.}

This is proved by approximating (in a certain sense) the Riemannian metric
by a triangulation, as follows.

\noindent {\bf Proposition 5.1.} {\sl For each $K$, $I > 0$, there is a constant
$c'(K,I)$ with the following property. If $M$ is a compact 3-manifold with a
Riemannian metric having injectivity radius at least $I$ and all sectional curvatures
in the interval $(-K,K)$, then $M$ has a triangulation $T$ such that
\item{1.} there is a $c'(K,I)$-Lipshitz homeomorphism $M_{PL} \rightarrow M$, where
$M_{PL}$ is the path metric on $M$ obtained by realising each tetrahedron of $T$
as a standard Euclidean simplex with side length $1$, and
\item{2.} the number of tetrahedra in $T$ is at most
$c'(K,I) \ V$, where $V$ is the volume of $M$.

}

\noindent {\sl Proof.} For each point $x$ in $M$, let ${\rm exp}_x \colon T_xM \rightarrow M$
be the exponential map. If we set $\epsilon > 0$ to be small enough (as a function of
$K$ and $I$), then the restriction of ${\rm exp}_x$ to the ball of radius $\epsilon$
about $0$ is injective for all $x \in M$. Furthermore, if $\epsilon$ is sufficiently small, then
$(1/\epsilon) B_{\epsilon}(x)$ is nearly isometric to the unit ball in
${\Bbb R}^3$, where $B_{\epsilon}(x)$ is the
ball of radius $\epsilon$ about $x$, and $(1/\epsilon) B_{\epsilon}(x)$ is obtained 
from $B_{\epsilon}(x)$ by rescaling by the factor $(1 / \epsilon)$.
Thus, if $x$ and $x'$ are within distance $\epsilon$ of each other,
then the set of points in $B_\epsilon(x)$ that are equidistant from
$x$ and $x'$ is a smooth surface, and when it is rescaled by the factor 
$(1/\epsilon)$, it is nearly isometric to a subset of the Euclidean plane.

Pick a maximal set of points $\{ x_1, \dots, x_n \}$ in $M$, such that no two of
these points are within $\epsilon/2$ of each other. Then, by maximality,
every point of $M$ lies at a distance of at most $\epsilon/2$ from at least one $x_i$.

Now consider the cut locus $L$ of these points. More precisely,
for each point $x$ in $M$, define $N(x)$ to be the number of
$\{ x_1, \dots, x_n \}$ that are closest to $x$. Then $L$
is the set of points $x$ where $N(x) \geq 2$.
By a small perturbation of the points $\{ x_1, \dots, x_n \}$, we may assume
that $L$ is a cell complex, where the set of points $x$ with $N(x) = k$ is 
the open $(4-k)$-cells.
A schematic picture of $L$, reduced by one dimension, is shown on the left in Figure 8.

\vskip 6pt
\centerline{
\epsfxsize=2.8in
\epsfbox{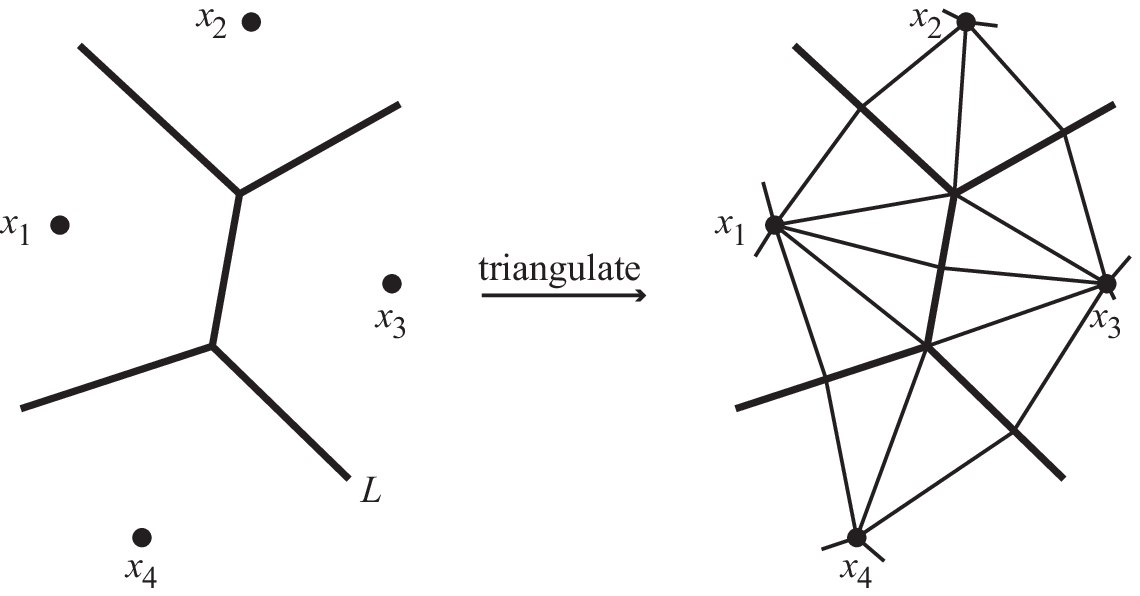}
}
\vskip 6pt
\centerline{Figure 8.}

Subdivide each 1-cell by adding in its midpoint.
Each 2-cell, rescaled by $(1/\epsilon)$, is nearly isometric to
a Euclidean polygon. In particular, if $\epsilon$ is sufficiently small,
it is star-shaped about some point in the interior. Subdivide
the 2-cell by adding a vertex at this point, and then coning off
the boundary. Now subdivide each 3-cell by coning off from the
point in $\{ x_1, \dots, x_n \}$ that it contains. The result is
the triangulation $T$.

We now wish to bound the number of tetrahedra in $T$. Each tetrahedron
has a vertex at one of $\{ x_1, \dots, x_n \}$. So our first step is to
find an upper bound for $n$. The points $x_1, \dots, x_n$ are more than
$\epsilon/2$ apart, and so the balls of radius $\epsilon/4$ about these
points are all disjoint. There is a lower bound $c_1(K,I)$ on the volume of each such
ball. Hence, $n \leq V/c_1(K,I)$, where $V$ is the volume of $M$.

So, our next task is to bound the number of tetrahedra incident to some $x_i$.
Each tetrahedron has a vertex in the interior of a 2-cell. This 2-cell
is equidistant between $x_i$ and some $x_j$. These are the closest points
in $\{ x_1, \dots, x_n \}$ to that vertex. But every point of $M$ is at most
$\epsilon/2$ from some point of $\{ x_1, \dots, x_n \}$. Hence, $x_j$ is
at most $\epsilon$ from $x_i$. By our assumption on the $\epsilon$,
the points $x_i$ and $x_j$ determine just one 2-cell.
The number of such $x_j$ is at most
$${\rm Volume}(B_{5\epsilon/4} (x_i))/ \min\{{\rm Volume}(B_{\epsilon/4}(x_j)) : 1 \leq j \leq n\},$$
because the $\epsilon/4$ balls about these $x_j$ are all disjoint, and
they fit into the $5\epsilon /4$ ball about $x_i$. Hence, the number
of such $x_j$ is at most some constant $c_2(K,I)$.

The remaining two vertices of any tetrahedron that is incident to $x_i$
lie in the boundary of the 2-cell, at either a corner or the midpoint of an edge. 
Again, the number of such vertices is at most some constant $c_3(K,I)$. Hence, putting this all together,
we deduce that the number of tetrahedra of $T$ is at most
$c'(K,I) \ V$, for some constant $c'(K,I)$.

We now construct the Lipschitz homeomorphism $M_{PL} \rightarrow M$.
Each simplex of $M_{PL}$ is already realised as a subset of $M$.
So, we only need to specify how the Euclidean simplex in $M_{PL}$
is mapped into that subset.
We start by sending the vertices of $M_{PL}$ to the corresponding points in $M$.
Then we consider the edges of the triangulation that lie in 1-cells of $L$.
We map these to $M$ so that they have constant speed. Then we consider faces
of the triangulation which lie in 2-cells of $L$. These have one vertex
in the interior of the 2-cell and two in the 1-skeleton. We view the
face as a cone on the former vertex and map this into $M$ in a way that
respects the cone structure. In other words, each point of the face lies on
a unique geodesic in the Euclidean metric on the face that runs from the coning
vertex to the opposite edge. We send this geodesic to the corresponding
curve in the 2-cell of $L$ that is a geodesic in the path metric on $L$.
We then extend the map over the remainder of each tetrahedron of $M_{PL}$
by viewing it as a cone with cone point being the vertex in the interior of a 3-cell. 
It is clear that this map
is $c'(K,I)$-Lipschitz for some constant $c'(K,I)$. Note that it need not
be bi-Lipschitz, because it is possible that some tetraheda in $M_{PL}$ may be
mapped to nearly flat tetrahedra in $M$. $\square$

\noindent {\sl Proof of Theorem 1.3.} Let $M$ be a solid torus with a Riemannian
metric as in the statement of the theorem. Let $T$ be the triangulation given by Proposition 5.1.
By Theorem 1.1, there is a pre-core curve $C$ that lies in the 2-skeleton of $T$,
that intersects the 1-skeleton in only finitely many points and that intersects the interior of
each face in at most 10 straight arcs. We make an arbitrarily small
perturbation of $C$ to make it into a core curve. 
So the length of $C$ in $M_{PL}$ is at most $40 |T|$, where $|T|$
is the number of tetrahedra. By (2) of Proposition 5.1, $|T|$ is at most $c'(K,I) V$.
By (1) of Proposition 5.1, there is a Lipschitz homeomorphism $M_{PL} \rightarrow M$ with Lipschitz constant
$c'(K,I)$. Thus, the length of the image of $C$ in $M$ is at most $40 (c'(K,I))^2 V$. $\square$

\vskip 18pt
\centerline{\caps 6. Examples}
\vskip 6pt

In this section, we investigate a family of triangulations of the solid torus.
Using them, we will prove the following.

\noindent {\bf Theorem 6.1.} {\sl There exists a family of triangulations
$T_i$ ($i \in {\Bbb N})$ of the solid torus $M$, with the following properties:
\item{1.} the number
of triangles and squares in any normal meridian disc is at least $c^{|T_i|}$, where
$|T_i|$ is the number of tetrahedra in $T_i$, and $c$ is the golden ratio
$(1 + \sqrt 5)/2$;
\item{2.} each pre-core curve in $\partial M$ intersects the 1-skeleton
of $T_i$ in at least $c^{|T_i| - 2}$ points;
\item{3.} there is a pre-core curve that lies in the 2-skeleton of $T_i$ and that intersects
the 1-skeleton in just one point. Moreover, when $i \geq 1$, this is in fact a core curve.

}

These examples exhibit the intrinsically exponential nature of normal surface theory.
They also show that it would be impossible to prove Theorem 1.1 by using only pre-core
curves in $\partial M$.

\noindent {\sl Proof.} Start with a fixed triangulation $T_0$ of the solid torus $M$. For the sake of being definite,
we use the triangulation with a single tetrahedron. (This is described in 
[4] for example.) This restricts to
a one-vertex triangulation of the boundary torus. The three edges 
on the boundary have slopes $(1,0)$, $(1,1)$ and $(2,1)$. Here, we
are using the standard basis for the first homology of the torus,
where $(0,1)$ is a meridian and $(1,0)$ is a longitude.

We will construct the triangulations $T_i$ recursively, with
each obtained from its predecessor by attaching a tetrahedron onto
its boundary, as shown in Figure 9. This has the effect of
performing an elementary move on the boundary triangulation.
This removes one of the three edges of
the boundary triangulation, so that the two triangles patch together to form
a square, and then inserts the other diagonal of this square.

\vskip 18pt
\centerline{
\epsfxsize=3in
\epsfbox{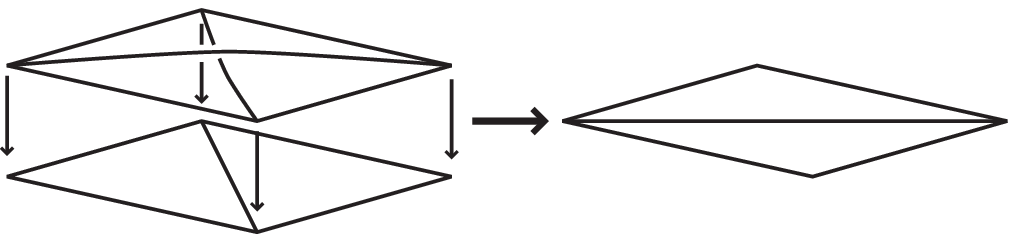}
}
\vskip 6pt
\centerline{Figure 9.}

\vskip 18pt
\centerline{
\epsfxsize=3in
\epsfbox{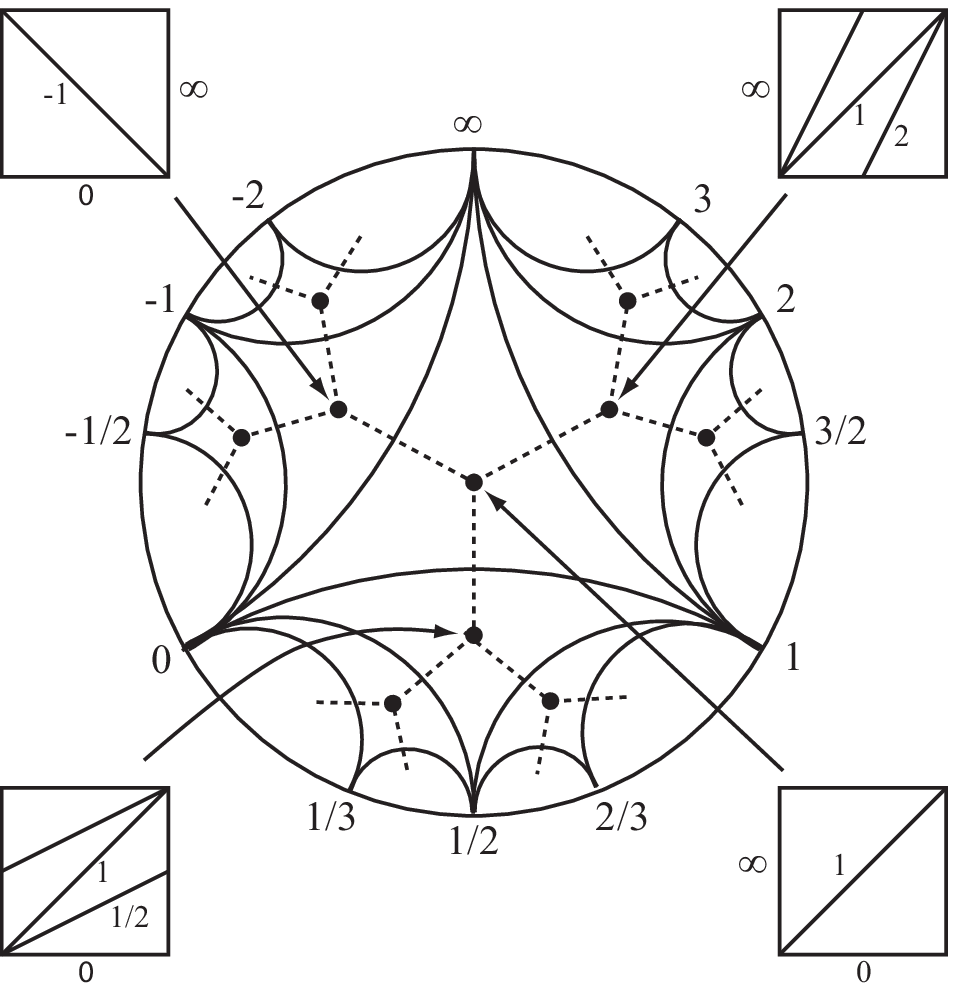}
}
\vskip 6pt
\centerline{Figure 10.}
\vfill\eject

It is well known that the set of one-vertex triangulations of a torus
forms the vertices of a tree. Two vertices are joined by an
edge if and only if the corresponding triangulations differ by
an elementary move. See Figure 10.
Our initial triangulation $T_0 \cap \partial M$ corresponds to one of these vertices
(the one just below the centre). To obtain $T_1$,
we perform the elementary move that removes $(1,0)$ and inserts
$(3,2)$. Then, to obtain $T_2$, we remove $(1,1)$ and insert $(5,3)$.
We repeat in this way, following  the path in this tree that turns left,
then turns right, then turns left, and so on.
This gives our sequence of triangulations $T_i$ of $M$.

We now compute the slopes of the 1-cells in
the boundary torus of $T_i$.
Let $s_i$ be the slope of the 1-cell that is removed when passing from
$T_i$ to $T_{i+1}$. We orient this slope so that it represents a
homology class with non-negative intersection number with the meridian.
Then each $s_i$ is, in standard homology co-ordinates, $(x_i, y_i)$.
So, $s_0 = (1,0)$, and $s_1 = (1,1)$. It follows
from our construction that, for each $i$, $s_{i+2} = s_i + s_{i+1}$.
In particular, the integers $y_i$ satisfy the Fibonacci relation
$y_{i+2} = y_{i} + y_{i+1}$, where $y_0 = 0$ and $y_1 = 1$.
Hence,
$$y_i = {1 \over \sqrt 5} \left ( {1 + \sqrt 5 \over 2} \right )^i 
- {1 \over \sqrt 5} \left ( {1 - \sqrt 5 \over 2} \right )^i.$$
The integers $x_i$ also satisfy the Fibonacci relation, and
since $x_0 = 1 = y_1$ and $x_1 = 1 = y_2$, we deduce that $x_i = y_{i+1}$
for all $i$.

We want to find a lower bound on the number of triangles and
squares of any normal meridian disc in $T_i$. Now, $T_i$ contains
an edge with slope $s_{i+2}$ and each normal triangle or square can intersect
this edge at most once. Thus, it suffices to find a lower bound on
the intersection number between $s_{i+2}$ and any meridian curve.
But this is $x_{i+2}$, which is at least $c^{i+1} = c^{|T_i|}$.
This proves (1).

Proving (2) is slightly more tricky because a curve on the boundary
of the solid torus with slope $(1,n)$, for any integer $n$, is a
pre-core curve. However, the intersection number between $(1,n)$
and $s_{i+2}$ is
$$|n x_{i+2} - y_{i+2}| = |n x_{i+2} - x_{i+1}| = x_{i+2} \ |n - (x_{i+1}/x_{i+2})|.$$
Now, as $i \rightarrow \infty$, $x_{i+1}/x_{i+2} \rightarrow c^{-1}$ which is
not an integer. So, there is a uniform lower bound to the difference
$|n - (x_{i+1}/x_{i+2})|$ provided $i$ is sufficiently large, and in fact one can verify that 
it is always at least $1/3$. So, the intersection number between
any pre-core curve on $\partial M$ and $s_{i+2}$ is at least
$x_{i+2}/3 \geq c^{i-1} = c^{|T_i|-2}$, as required. 

Finally, note that $T_0$ contains a pre-core curve that lies
in the 2-skeleton and that intersects the 1-skeleton in just one point.
This point lies in the interior of the edge on $\partial M$ with slope $(1,0)$. Thus,
this curve becomes a core curve in $T_1$ and all subsequent triangulations. $\square$

\vskip 18pt
\centerline{\caps References}
\vskip 6pt

\item{1.} {\caps W. Haken,} {\sl Theorie der Normalfl\"achen.}
Acta Math. 105 (1961) 245--375. 

\item{2.} {\caps J. Hass, J. Lagarias}, {\sl The number of Reidemeister moves needed
for unknotting}, J. Amer. Math. Soc. 14 (2001) 399--428.

\item{3.} {\caps W. Jaco, U. Oertel}, {\sl An algorithm to decide if a $3$-manifold 
is a Haken manifold,} Topology 23 (1984) 195--209. 

\item{4.} {\caps W. Jaco, J. H. Rubinstein, S. Tillmann},
{\sl Minimal triangulations for an infinite family of lens spaces},
J. Topol. 2 (2009) 157--180. 

\item{5.} {\caps M. Lackenby}, {\sl Word hyperbolic Dehn surgery}, Invent. Math. 140 (2000) 243--282.

\item{6.} {\caps M. Lackenby}, {\sl The crossing number of composite knots}, J. Topol.
 2 (2009) 747--768.

\item{7.} {\caps M. Lackenby}, {\sl The crossing number of satellite knots}, To appear.

\item{8.} {\caps S. Matveev,} {\sl Algorithmic topology and classification of $3$-manifolds},
Algorithms and Computation in Mathematics, Volume 9, Springer (2003). 

\vskip 18pt
\+ Mathematical Institute, University of Oxford, \cr
\+ 24-29 St Giles', Oxford OX1 3LB, United Kingdom. \cr

\end